\documentclass[12pt]{article}
\usepackage{latexsym}
\usepackage{amssymb}
\usepackage{amsmath}
\usepackage{amsfonts}
\newcommand{\N}{\mathbb{N}}

\newtheorem{theorem}{Theorem}

\title{Applications of hypercomplex automorphic forms in Yang-Mills gauge theories}
\author{
Rolf S\"oren~Krau{\ss}har
\thanks{Erziehungswissenschaftliche Fakult\"at, Lehrgebiet Mathematik und ihre Didaktik, Universit\"at Erfurt, 
D-99089 Erfurt, Germany.
E-mail: {\tt soeren.krausshar@uni-erfurt.de} } \and J\"urgen
Tolksdorf
\thanks{Max-Planck Institute for Mathematics
in the Sciences, Inselstra{\ss}e 22, D-04103 Leipzig, Germany.
E-mail: {\tt juergen.tolksdorf@mis.mpg.de}}}
\begin{document}
\maketitle
\begin{abstract}
In this paper we show how hypercomplex function theoretical
objects can be used to construct explicitly self-dual
$SU(2)$-Yang-Mills instanton solutions on certain classes of
conformally flat $4$-manifolds. We use a hypercomplex
argument principle to establish a natural link between the
fundamental solutions of $D \Delta f = 0$ and the second Chern
class of the $SU(2)$ principal bundles over these manifolds. The
considered base manifolds of the bundles are not simply-connected,
in general. Actually, this paper summarizes an extension of the
corresponding results of G\"ursey and Tze on a hyper-complex
analytical description of $SU(2)$ instantons. Furthermore, it provides an application of the recently introduced new 
classes of hypercomplex-analytic automorphic forms.
\end{abstract}



\section{Basic concepts of quaternionic analyticity in Euclidean Space}

In this section we summarize some basics on quaternions and their function theory.
For more details we refer the reader, for instance, to \cite{BDS,DSS}. Basic
applications of this framework to $SU(2)$-Yang Mills gauge theories may be found in 
the seminal works by  G\"ursey and Tze \cite{GT1,GT2}.
\par\medskip\par
Let $\mathbb{H} \stackrel{\sim}{=}\mathbb{R} \oplus \mathbb{R}^3$
be the set of real quaternions. Each quaternion $ a \in \mathbb{H}$ is
written as $a = \sum_{\mu = 0}^3 a_{\mu} e_{\mu} = a_0 e_0 + \sum_{i=1}^3 a_i e_i \equiv
Sc(a) + Vec(a),\; a_\mu \in \mathbb{R}$. Here, $Sc(a) = a_0$ is called the scalar or real part 
of the quaternion $a$; $Vec(a) = \sum_{i=1}^3 a_i e_i$ is called the vector or imaginary part of $a$.

The quaternionic
multiplication is defined by $e_i^2=-1$ for $i=1,2,3,$ $e_0^2=1$
and $e_1 e_2 = e_3 = -e_2 e_1$, $e_2 e_3 = e_1 = -e_3 e_2$, $e_3
e_1 = e_2 = -e_1 e_3,$ as well as $e_0e_i = e_i e_0$ for
$i=1,2,3$.

The quaternionic conjugate is given by  $ \overline{a} := Sc(a) - Vec(a).$ 
Accordingly, the norm is defined as $N(a)\equiv |a|^2 := a \overline{a} =
\overline{a} a = \sum_{\mu=0}^3 a_\mu^2$. It coincides with the Euclidean norm of $\mathbb{R}^4$,
such that as Euclidean spaces $\mathbb{H}\simeq\mathbb{R}^{4,0}$. 
\par\medskip\par
Next, we call in mind two different concepts of quaternionic analyticity. For more details on quaternionic analysis we refer, for instance, to \cite{BDS,DSS,GS1,GS2}. For new connections between the different concepts of analyticity in hypercomplex spaces we also refer to the recent works \cite{CGK,CGS}. 
\par\medskip\par
Let $U \subset \mathbb{H}\simeq\mathbb{R}^{4,0}$ be an  open subset. A quaternion-valued function 
$f:U \rightarrow\mathbb{H}$, which is real differentiable in each real component, is called 
{\it left monogenic} (resp. {\it right monogenic}) on $U$ if $D f = 0$ (resp. $fD =0$). Here, the first
order differential operator 
$$D := \frac{\partial}{\partial x_0} + 
e_1 \frac{\partial}{\partial x_1} + e_2 \frac{\partial}{\partial x_2} + 
e_3 \frac{\partial}{\partial x_0} $$ 
is the quaternionic analogue of the Cauchy-Riemann operator 
${\cal{D}} = \frac{\partial}{\partial x} + i \frac{\partial}{\partial y}$ of complex analysis.

The operator $D$ often abbreviated as   
$\partial_0 + e_1
\partial_1  + e_2 \partial_2 + e_3 \partial_3$. In what follows we focus on the class of
left monogenic functions and will call them simply monogenic for short.

\par\medskip\par
An important property of the quaternionic Cauchy-Riemann operator is, that is factorizes the Euclidean Laplacian
$$
\Delta = \sum\limits_{i=0}^3 \frac{\partial^2}{\partial x_i^2}
$$
viz $\Delta = D \overline{D} = \Big(\partial_0 + e_1
\partial_1  + e_2 \partial_2 + e_3 \partial_3\Big)\Big(\overline{\partial_0 + e_1
\partial_1  + e_2 \partial_2 + e_3 \partial_3}\Big)$. 
\par\medskip\par

Each quaternion $a \in \mathbb{H}$ is known to have a polar decomposition: $a = a_0 + \omega r$, 
with $r > 0$ and $\omega$ being a unit vector in $\mathbb{R}^3\subset \mathbb{H}$. Since 
$\omega^2 = -1$, the vector $\omega$ defines a complex structure on $\mathbb{R}^4$.

Let $F: \mathbb{C} \rightarrow \mathbb{C}, F(x + iy) = u(x,y) + i v(x,y)$ be a complex analytic function. 
It follows that 
$$G: \mathbb{H} \rightarrow \mathbb{H},\quad G(a_0 + \omega r) = u(a_0,r) + \omega v(a_0,r)$$ 
satisfies the linear third order equation $D \Delta G = 0$, see for instance \cite{Fu32}. The function $G$ 
is called (left) {\it Fueter-holomorphic}, cf. for instance \cite{Fu32,Sce,DSS,LR1,Sommen}. 
Note that the function $f:= \Delta G$ is left monogenic. This was proved already in the papers cited above.  

\par\medskip\par
{\bf Remark:} Since $\Delta$ is a scalar operator, one has $\Delta D  = D \Delta$ and therefore $\Delta D f = 0$ if and only if $D \Delta f = 0$, whenever $f$ is at least three times continuously differentiable.
\par\medskip\par

All integer powers of $z=x_{\mu} e_{\mu} = x_0 + x_1e_1 + x_2 e_2 + x_3 e_3 \in\mathbb{H}$
have the property that they are Fueter-holomorphic, cf. for instance \cite{LR1}. However, none
of these are left or right monogenic. As suggested in \cite{BDS} and elsewhere, in the monogenic
context the negative power function $\frac{1}{z}$ is replaced by the fundamental
solution of $D$ having the form 
\begin{equation}\label{q0}
q_{\bf 0}(z) =
\frac{\overline{z}}{|z|^4} = \frac{x_0 - x_1 e_1 + x_2 e_2 + x_3 e_3}{(x_0^2+x_1^2+x_2^2+x_3^2)^2}= - \frac{1}{4} \Delta z^{-1}.
\end{equation}
The other negative powers $z^{-m}$ with $m \in \mathbb{N}^{\ge 2}$ are replaced by the partial derivatives of $q_{\bf 0}(z) $. 
The partial derivatives of $q_{\bf 0}(z)$ will be denoted by 
$$
q_{\bf m}(z) := q_{m_1,m_2,m_3}(z) := \frac{\partial^{m_1+m_2+m_3}}{\partial x_1^{m_1}\partial x_2^{m_2}\partial x_3^{m_3} }q_{\bf 0}(z),
$$
where ${\bf m} := (m_1,m_2,m_3)$ is a multi index from $\mathbb{N}_0^3$. 

Consequently, all quaternionic M\"obius
transformations are Fueter-holomorphic, but none of them is left
monogenic. For convenience we recall that M\"obius transformations
can be expressed in the quaternionic language in the form
$$M(z)=(az+b)(cz+d)^{-1},$$ 
where $a,b,c,d$ are arbitrary quaternions satisfying the determinant condition (see for instance 
\cite{Ahl86,Zoll}):
$$
|b-ac^{-1}d||c| \neq 0,\;\;\mbox{if}\;\; c \neq 0 \;\;
\mbox{and}\;\; |ad| \neq 0 \;\;\mbox{if}\;\; c = 0.
$$
Due to Liouville's theorem (see for example \cite{Bl}), these are exactly the conformal maps 
in $\mathbb{H}$. 
\par\medskip\par
Both
operators $D$ and $D\Delta$ are conformally invariant up to an
automorphy factor. According to \cite{Fu48/49,KraHabil,Ry85,Ry93}, suppose that $D_y f(y) = 0$
resp. $D_y \Delta_y f(y) =0$ and $y=(az+b)(cz+d)^{-1}$ is a M\"obius transformation. Then 
$$q_{\bf 0}(cz+d)
f((az+b)(cz+d)^{-1}),$$ resp. 
$$(cz+d)^{-1} f((az+b)(cz+d)^{-1})$$
lies again in an appropriate domain of $\mathbb{H}$ in Ker $D_z$
resp. Ker $D_z \Delta_z$.

\par\medskip\par

Similar to complex function theory, as proved in \cite{BDS} and elsewhere, quaternion-valued functions in
${\rm Ker}(D)$ satisfy a Cauchy integral formula of the following form
$$
f(z) = \frac{1}{8 \pi^2}\int_{\partial K} q_{\bf 0}(z-w)
d\sigma(w) f(w),
$$
where $K$ is some compact subset with strongly Lipschitz boundary.
Furthermore, $K$ is contained in an open set $U$, such that $f$ satisfies
$Df=0$. Here, and in all that follows
$$
d\sigma(w) = dw_1 \wedge dw_2 \wedge dw_3 - dw_0 \wedge dw_2
\wedge dw_3 e_1 + dw_0 \wedge dw_1 \wedge dw_3 e_2 - dw_1 \wedge
dw_2 \wedge dw_3 e_3
$$
denotes the oriented surface $3$-form. 
Notice that the expression $q_{\bf 0}(z-w)$ replaces the complex function $\frac{1}{z-w}$ in the 
ordinary Cauchy integral formula. In the quaternionic version the factor $8 \pi^2$ replaces the factor 
$2 \pi i$ from complex analysis.  
\par\medskip\par
From this generalized Cauchy integral formula one may establish a Green's integral
formula for functions in ${\rm Ker}(D \Delta)$, see
\cite{KraHabil}. We also need to consider the more general
topological version of Cauchy's integral formula, cf.
\cite{KraHabil}
\begin{theorem}
Let $U \subset \mathbb{H}$ be open. Let $\Gamma$ be a
$3$-dimensional null-homologous cycle in $U$. Suppose that $z \in U
\backslash \Gamma$ and that $f:U \rightarrow \mathbb{H}$ is left
monogenic. Then $$\frac{1}{8 \pi^2} \int\limits_{\Gamma} q_{\bf
0}(z - \zeta) d\sigma(\zeta)f(\zeta)= w_{\Gamma}(z)f(z),$$ where
$w_{\Gamma}(z)$ stands for the winding number of $\Gamma$ with
respect to $z$.
\end{theorem}
The winding number counts how often $\Gamma$ wraps around the point $z$.  
This generalized topological Cauchy formula provides us with a quaternionic 
analogue of the usual residue theorem in the form
\begin{equation}\label{res}
\int\limits_{\Gamma} d\sigma(z) f(z) = 8 \pi^2 \sum\limits_{i=1}^s
{\rm res}(f;\beta_i),
\end{equation}
where $f$ is supposed to be left monogenic in ${\rm int}\Gamma$ except
of in a finite number of isolated points denoted by $\beta_i$. 
\par\medskip\par
{\bf Remark}: In complete analogy to complex function theory, the residue of a function $f$ at a point $\beta_i$ is nothing else 
than the first Laurent coefficient of the singular part of the Laurent series expansion, i.e. the coefficient that is 
associated with $q_{\bf 0}$ in the expansion given below in (\ref{laurent}). Following for instance \cite{BDS,KraHabil}, the Laurent series of a function 
that is left monogenic in the open pointed ball $B(\beta_i,\varepsilon)\backslash \{\beta_i\}$ of radius $\varepsilon$ and 
center $\beta_i$  
is of the form
\begin{equation}\label{laurent}
f(z) = \sum\limits_{{\bf m} \in \mathbb{N}_0^3} q_{\bf m}(z-\beta_i) b_{\bf m} + \sum\limits_{{\bf m} \in \mathbb{N}_0^3} V_{\bf m}(z-\beta_i) a_{\bf m}
\end{equation}
where $V_{\bf m}$ are the so-called Fueter polynomials, $q_{\bf m}$ the partial derivatives of the Cauchy kernel 
and $a_{\bf m}$ and $b_{\bf m}$ quaternionic coefficients. Hence, the residue of $f$ at $\beta_i$ 
equals  ${\rm res}(f;\beta_i)=b_{\bf 0}$.    

In the special case where the function $f$ is left monogenic at each point of the open ball $B(\beta_i,\varepsilon)$, this series 
expansion simplifies to the Taylor type series expansion of the form 
$$
f(z) = \sum\limits_{{\bf m} \in \mathbb{N}_0^3} V_{\bf m}(z-\beta_i) a_{\bf m}.
$$
Consequently, this residue formula then simplifies to 
$$
\int\limits_{\Gamma} d\sigma(z) f(z) = 0,
$$
and we have a generalization of Cauchy's theorem in a ball, which can be extended to open star-like domains by applying standard 
arguments from the literature. 
\par\medskip\par
As a further consequence of the generalized Cauchy integral formula combined with the generalized Cauchy integral theorem, the following generalization of the complex argument principle was proved
a few years ago to hold for left monogenic functions with isolated $a$-points, cf. \cite{HeKra,KraHabil}. 

For convenience we also recall that $c$ is an isolated $a$-point of $f$ if $f(c)=a$ and if additionally there exists a sufficiently small neighborhood $V$  around $c$ where $f(x) \neq a$ for all $x \in V \backslash\{a\}$.
\begin{theorem}({\rm cf. \cite{HeKra,KraHabil}}).

Let $G\subset \mathbb{H}$ be a domain. Suppose that $f:
G\rightarrow \mathbb{H}$ is left monogenic in $G$ and that $c\in
G$ is an isolated zero point of $f$.   
 
Next, let $\varepsilon > 0$ so that
$\overline{B}(c,\varepsilon)\subseteq G$ and
$f|_{\overline{B}(c,\varepsilon)\backslash\{c\}}\neq 0$. Then
$$
{\rm ord}(f;c) =\frac{1}{8 \pi^2}\int\limits_{\partial
B(c,\varepsilon)} q_{\bf 0}(f(z)) \cdot \Big[(Jf)^{ad}(z)\Big] *
\Big[d\sigma(z)\Big]
$$
\end{theorem}
Here $``*"$ denotes the matrix multiplication and $``\cdot"$ the
quaternionic multiplication. The vector
$[(Jf)^{ad}(z')]*[d\sigma'(z')]\in\mathbb{R}^4$ then is re-interpreted as
a quaternion. $(Jf)^{ad}$ stands for the adjunct matrix of the Jacobian of $f$. 

These tools will now be used to study self-dual $SU(2)$ Yang-Mills
instantons on specific classes of conformally flat $4$-manifolds. These new tools together with the new class of 
automorphic forms developed in \cite{KraHabil,BCKR,CGK} allow to round off some classical studies 
by G\"ursey and Tze on self-duality in $SU(2)$-Yang-Mills gauge theories and to review their results from the viewpoint of a new mathematical fundamental theory.

\section{Self-duality of Yang-Mills instantons and quaternionic analyticity}

\subsection{Self-duality in $SU(2)$-Yang-Mills gauge theory}

In \cite{GT1,GT2}, F. G\"ursey and Tze constructed an explicit relation between the self-duality condition 
for $SU(2)$-Yang-Mills instantons on $\mathbb{R}^4$ and Fueter holomorphic functions. To get started, 
we briefly summarize some of their results.
\par\medskip\par
For this we consider the quaternionic Hopf-bundle $Sp(1)
\hookrightarrow S^7 \twoheadrightarrow S^4$ with typical fiber $S^3\simeq SU(2)\subset\mathbb{H}.$ 
The base manifold $S^4\simeq\mathbb{H}{\rm P}_{\!1}$ is regarded as being the one-point
compactification of $\mathbb{R}^{4,0}\simeq\mathbb{H}.$ Moreover, the total space of the quaternionic Hopf-bundle $S^7$ is considered as being a sub-manifold of $\mathbb{H}\times\mathbb{H}$. In other words, the quaternionic Hopf-bundle is considered as being a natural sub-bundle of the trivial quaternionic line bundle $\mathbb{R}^4\times\mathbb{H}\twoheadrightarrow\mathbb{R}^4.$ On the latter, respectively, the gauge group ${\cal G}$ and the affine set of connections ${\cal A}$ can be identified with
$C^{\infty}(\mathbb{R}^4,SU(2))$ and $C^{\infty}(\mathbb{R}^4,
\Lambda^{\!1}\mathbb{R}^{4}\otimes su(2))$. Here, $su(2)$ denotes the Lie algebra of $SU(2).$ 

The field strength $F_{\!A}\in C^{\infty}(\mathbb{R}^4,\Lambda^{\!2}\mathbb{R}^{4}\otimes su(2))$ of
$A\in C^{\infty}(\mathbb{R}^4,\Lambda^{\!1}\mathbb{R}^{4}\otimes su(2))$ is defined by 
$F_{\!A}:= dA + A\wedge A$. As a consequence, the field strength satisfies the Bianchi-Identity: 
$dF_{\!A}=-[A,F_{\!A}].$ The gauge group ${\cal G}$ naturally acts on ${\cal A}$ from the right via
$$
{\cal{A}} \times {\cal{G}} \rightarrow {\cal{A}},\; (A,g) \mapsto
A^g := g^{-1} A g + g^{-1} dg.
$$
It follows that $F_{\!\!A}^g\equiv F_{\!A^g} = g^{-1}F_{\!A}\,g$. The Yang-Mills functional
$$
S_{\mbox{\tiny YM}}: {\cal{A}} \rightarrow \mathbb{R},\; A \mapsto \int_{S^4}
tr(F_{\!A} \wedge\ast F_{\!A})
$$
is thus invariant with respect to the right-action of ${\cal G}$ on ${\cal A}$. Here, ``$tr$'' refers to the ordinary matrix trace and ``$\ast$'' denotes the Hodge map with respect to the orientation and the 
Fubini-Study metric on $\mathbb{H}{\rm P}_1.$.

The critical points of the Yang-Mills functional fulfill the {\it Yang-Mills equation} 
$d\ast F_{\!A} = -[A,\ast F_{\!A}].$ In the definition of the Yang-Mills functional the Yang-Mills field 
strength $F_{\!A}$ is identified with its pull-back with respect to the stereographic
projection $S^4\twoheadrightarrow\mathbb{R}^4.$ Of particular interest are the solutions 
$F_{\!A}$ satisfying $F_{\!A}\pm *F_{\!A} = 0$. These are the anti-self dual (resp. self-dual) 
{\it instanton solutions}. Geometrically $A\in{\cal A}$ is interpreted as a gauge potential of an 
$SU(2)$-connection on the Hopf bundle.

\par\medskip\par

In \cite{GT1,GT2} it has been shown that $$a_{\mu} = \frac{1}{2}
\frac{t {\stackrel{\leftrightarrow}{\partial}}_{\mu}t}{1+t^{\sharp}t}, \quad t \in \mathbb{H},$$ transforms 
under the group ${\cal C}^\infty(\mathbb{R}^4,Sp(1))$ as: 
$a_{\mu} \rightarrow \overline{r} a_{\mu} r + \overline{r}\partial_{\mu} r$. This means that 
$A:=dx^\mu\otimes a_{\mu}$ transforms like a gauge potential under $Sp(1) = SU(2)$.

\par\medskip\par

Hence, one may define an $SU(2)$ Yang-Mills field
$F_{\!A}:=\frac{1}{2}dx^\mu\wedge dx^\nu\otimes f_{\mu\nu}$ with 
$f_{\mu\nu} := \partial_{\mu} a_{\nu} - \partial_{\nu} a_{\mu} + [a_{\mu},a_{\nu}]$. Because of the 
Bianchi identity, the Yang-Mills equation are automatically satisfied if $F_{\!A}$ is self-dual or anti-self-dual.

Next, one may introduce the Fueter holomorphic function
$${\cal F}(z) = (z+a)^{-1} + \overline{\alpha}(z+\beta)^{-1}\alpha,$$
where $a,\beta,\alpha$ are some constant quaternions. From
\cite{GT1,GT2} one infers that
$$
a_{\mu} = \frac{1}{2}Vec[e_{\mu}(\Delta{\cal F})(D{\cal F})^{-1}].
$$
Let $\Phi_{\mu\nu} = f_{\mu\nu} + \tilde{f}_{\mu\nu}$ denote the
self-dual part of $f_{\mu\nu}$ and $\Phi'_{\mu\nu} = f_{\mu\nu} -
\tilde{f}_{\mu\nu}$ the anti-self-dual part of $f_{\mu\nu}$. Here,
$\tilde{f}_{\mu\nu}$ are the components of the Hodge dual $\ast
F_{\!A}$ of $F_{\!A}.$ It follows that
$$\Phi'_{\mu\nu} = - \frac{1}{2} e'_{\mu\nu}(D{\cal F})^{-1} (D\Delta {\cal F})\quad \mbox{and}\quad
\Phi_{\mu\nu} = - \frac{1}{2}(D{\cal F})(D e_{\mu\nu}
\overline{D}) (D{\cal F})^{-1}.$$ The self-duality condition
$\Phi'_{\mu\nu} \equiv 0$ thus is satisfied if $D \Delta {\cal{F}} = 0$. The
self-duality condition may thus be rephrased in terms of Fueter
holomorphy for the special ansatz proposed by G\"ursey and Tze. Note that $(z+a)^{-1}$ actually is the Green's kernel
of $D \Delta$ up to a constant.

\subsection{The relation between the Chern-Pontryagin index and the
quaternionic winding number}

In \cite{GT1,GT2} it has been shown that the second Chern class is
proportional to 
$$ \Pi = -Sc(f_{\mu\nu}\tilde{f}_{\mu\nu})  = \frac{1}{4} \Delta[ \Delta \ln (D{\cal F})].$$ 
The density $\ast\Pi$ coincides with the Yang-Mills Lagrangian density in the case of self-dual instantons. 
\par\medskip\par
Next let $V \subset \mathbb{H}$ be a four-dimensional domain having a three dimensional Lipschitz boundary $\partial V$ that 
is homeomorphic to the three-dimensional unit sphere $S^3 = \{z \in \mathbb{H} \mid |z|=1\}$. 
\par\medskip\par
Following \cite{GT1,GT2}, the second Chern number of the $SU(2)$ principal
bundle reads
\begin{eqnarray*}
c_2 &=& \frac{1}{8 \pi^2} \int_V \Pi dx_0 dx_1 dx_2 dx_3\\ 
&=& \frac{1}{8 \pi^2} \int_V\Big( \frac{1}{4} \Delta[ \Delta \ln (D{\cal F})]\Big) dx_0 dx_1 dx_2 dx_3 \\
&=&\frac{1}{32 \pi^2} \int_{\partial V} d\sigma(z) [\Delta
\overline{D} \ln (D{\cal F})].
\end{eqnarray*}
In the previous line we have applied the classical Gauss-Ostrogradskii-Green formula. 
 
The second Chern number is a topological invariant. As a consequence, the
topological deformation of the three-dimensional contour $\partial V$ reduces to the sum of
small hyper-spheres (or cycles) $\partial B(\beta_i,\varepsilon)$ surrounding each pole (denoted 
by $\beta_i$) of the integrand $[\Delta\overline{D} \ln (D{\cal F})]$. In the 
function theoretic language this topological invariance is rephrased as the application of the Cauchy integral 
theorem. This allows one to reduce the integration over sufficiently small spheres of radius $\varepsilon$ surrounding the poles, instead of extending the integration over the full boundary contour $\partial V$. In general, the latter is very difficult to parametrize. To be more precise, we have 
\begin{eqnarray*}
c_2 &=&
\frac{1}{32 \pi^2} \int_{\partial V} d\sigma(z) [\Delta
\overline{D} \ln (D{\cal F})] \\
&=& \frac{1}{32 \pi^2} \int_{\cup_{i=0}^n \partial B(\beta_i,\varepsilon)} d\sigma(z) \Delta
\Big(\sum_{i=0}^n \Big(\beta_i-z\Big)^{-1}\Big) \\&=& \frac{1}{8 \pi^2}
\int_{\cup_{i=0}^n\partial B(\beta_i,\varepsilon)} d\sigma(z) \sum_{i=0}^n q_{\bf 0}(z-\beta_i). 
\end{eqnarray*}
In the previous equality we exploited the relation (\ref{q0}), namely that 
$$
q_{\bf 0}(z) = - \frac{1}{4} \Delta\Big( z^{-1}\Big).
$$
Applying in the next step a
suitable M\"obius transformation, one can remove one of these $n+1$ poles from
$\overline{D} \ln D{\cal F}$, namely by sending for instance $\beta_0$ to $\infty$, where again $\infty$ 
is interpreted in the sense of the one-point compactification of $\mathbb{H}$. The other quaternions 
$\beta_1,\ldots,\beta_n$ are then sent to other quaternions, say $\gamma_1,\ldots,\gamma_n$. 

By the application of the quaternionic residue theorem, see (\ref{res}), one finally obtains 
$$ 
c_2= \frac{1}{8
\pi^2} \int_{\cup_{i=1}^n\partial B(\gamma_i,\varepsilon)} d\sigma(z) \sum_{i=1}^n q_{\bf
0}(z-\gamma_i) = \sum_{i=1}^n {\rm res} (q_{\bf 0};\gamma_i) = n.$$

\section{Quaternionic-analytic automorphic forms and instanton solutions on conformally flat manifolds}

On the $4$-sphere and the projective space $\mathbb{C}P_2$
self-dual solutions of a $SU(2)$-Yang-Mills instantons are well studied. In this section we discuss how hypercomplex analysis may be used to explicitly construct self-dual instanton solutions on various other classes of conformally flat spin manifolds. In recent works \cite{KraRyan2,BCKR} we have been 
able to construct monogenic spinor bundles over large classes of examples of such manifolds in terms of hypercomplex-analytic Eisenstein- and Poincar\'e series that have been introduced in \cite{KraHabil} in a different context. Our main goal in this paper is to develop an explicit link between the fundamental solutions of $D\Delta$ on some manifolds and the Chern number of $SU(2)$ principal bundles over these manifolds.
\par\medskip\par
Conformally flat manifolds are Riemannian manifolds with vanishing Weyl tensor.
Equivalently, conformally flat manifolds are manifolds with an atlas whose transition functions are M\"obius transformations. Remember that 
M\"obius transformations coincide exactly with the  conformal maps
in $\mathbb{H}$. Form this viewpoint (cf. \cite{KraRyan2})
conformally flat manifolds can be regarded as
higher dimensional analogues of holomorphic Riemann
surfaces. Following for example the classical paper \cite{Kuiper}, a large class of examples
can be constructed by factoring out a domain $U \subseteq
\mathbb{H}$ by a discrete group of the M\"obius group. Then
${\cal{M}} = U/\Gamma$ is a conformally flat manifold. The
simplest examples are of course $\mathbb{H}
\stackrel{\sim}{=}\mathbb{R}^4$ and $S^4$. Further important examples
are given by
\begin{itemize}
\item
$C_p = \mathbb{H}/{\cal{T}}(\mathbb{Z}^p)$, where
${\cal{T}}(\mathbb{Z}^p) = \Big\langle \left(\begin{array}{cc} 1 &
e_0 \\ 0 & 1 \end{array}\right),\ldots, \left(\begin{array}{cc} 1
& e_{p-1} \\ 0 & 1 \end{array}\right) \Big\rangle\;\; 1 \le p \le
4$ is the discrete translation group associated to the lattice
$\mathbb{Z}^p$.  If $p \le 3$, then $C_p$ are $p$-cylinders and
for $p=4$: $T_4 = \mathbb{H}/{\cal{T}}(\mathbb{Z}^4)$ is the flat
$4$-torus.
\item
Let $U=\mathbb{H} \backslash \{0\}$ and $\Gamma=\{m^k,k \in
\mathbb{Z}\}$ which is a discrete dilatation group. Then,  $U
/\Gamma \stackrel{\sim}{=} S^3 \times S^1$ is the Hopf manifold
\item
$U=H^+(\mathbb{H}) = \{z \in \mathbb{H}, Sc(z) > 0\}$.

Let $\Gamma_p = \Big\langle \left(\begin{array}{cc} 1 & e_1 \\ 0 &
1\end{array} \right),\ldots,\left(\begin{array}{cc} 1 & e_p \\ 0 &
1\end{array}\right),\left(\begin{array}{cc} 0 & -1
\\ 1 & 0
\end{array} \right)\Big\rangle$ ($p < 4$) denote the special hypercomplex
modular group as discussed in \cite{KraHabil}. For $N \in \mathbb{N}$
consider
$$
\Gamma_p[N] = \Big\{\left(\begin{array}{cc} a & b \\ c & d
\end{array}\right) \in \Gamma_p,\;|\; a-1,b,c,d-1 \in N
\mathbb{Z}^4 \Big\}\,.
$$
The associated class of conformally flat manifolds contains
$k$-handled cylinders $(p < 3$) and $k$-handled tori ($p=3$).

Instead of the standard orthonormal lattice we can also consider some other choices of lattices, for instance the set of reduced Hurwitz 
quaternions where the basis $\{e_1,e_2,e_3\}$ is replaced by $\{e_1,e_2,\frac{e_1+e_2+e_3}{2}\}$. This lattice plays a crucial 
rule in analytic number theory, see for instance 
\cite{Krieg88}.  
This lattice is more dense than the standard one, and 
the fundamental domain of the associated hypercomplex modular group hence has a smaller volume.   
\end{itemize}
In all these cases $U \subset \mathbb{H}$ is the universal
covering space of a manifold ${\cal{M}} = U / \Gamma$. Hence, there
exists a unique projection map $p: U \rightarrow {\cal{M}}, x
\mapsto x\;$(mod $\Gamma$). Let $x':=p(x)$. If $A \subseteq
\mathbb{H}$ is open, then $A'=p(A)$ is open on ${\cal{M}}$.
Provided $f$ is $\Gamma$-invariant on $U$, then $f'=p(f)$ represents a well-defined spinor-valued section on ${\cal{M}}$. Notice that 
different decompositions of the period lattice may give rise to different choices of 
spinor bundles over ${\cal{M}}$, but for details about this topic we refer the reader to \cite{KraRyan2,BCKR} where this 
issue has been adressed extensively.

Furthermore,
$D'=p(D)$ and $\Delta'=p(\Delta)$ yields the hyper-complex
analogue of the Cauchy-Riemann operator and the Laplace operator on ${\cal{M}}$.
\par\medskip\par
In the spirit of \cite{GT1,GT2}, 
self-dual $SU(2)$ instantons on these manifolds are induced by
$\Gamma$-``periodizations'' of the self-dual instantons that we
have on Euclidean space. Consider the Yang-Mills connection
$$a_{\mu} = \frac{1}{2} Vec\{e_{\mu}(\Delta{\cal F})(D{\cal
F})^{-1}\}$$ associated with superpositions of
$\Gamma$-``periodic'' actions on the basic Fueter holomorphic function $z^{-1}$.
Their partial derivatives then provide  us in particular with
further self-dual instanton solutions on the associated manifolds. The key
question that arises in this context immediately is to ask for the
fundamental solutions  in terms of explicit formulas of $D$ (resp.
of $D \Delta$) on these manifolds. In the mainly mathematical
oriented works~\cite{KraRyan2,BCKR} an explicit answer to this question has been
given. These results will now be interpreted in terms of self-dual $SU(2)$-instantons on 
these manifolds and Chern numbers.

\subsection{Conformally flat cylinders and tori}

Following \cite{KraRyan2}, the fundamental solution of $D$ (resp
$D \Delta$) on $\mathbb{H}/\mathbb{Z}^p$ ($1 \le p \le 3$) are the
generalized monogenic (resp. Fueter holomorphic) cotangent
functions. The simplest example $p=1$ is $$
\cot^{(1)}_{D\Delta}(z) = (\beta z)^{-1} +
\sum_{n=1}^{+\infty}\Big[(\frac{\beta}{2\pi} z -n)^{-1} +
(\frac{\beta}{2 \pi}z+n)^{-1}\Big].$$ In the spirit of \cite{GT2},
this function is physically interpreted as infinite superpositions
of 't Hooft instantons along the $x_0$-axis placed at periodic
intervals. Notice that in the mathematical context, 't Hooft instantons are hypermonogenic functions \cite{L,Qber}. 
Further special instanton solutions on this cylinder with values in the trivial spinor bundle 
are obtained by considering derivatives:
$$
\partial_{x_0}^{m_0} \partial_{x_1}^{m_1} \partial_{x_2}^{m_2}\partial_{x_3}^{m_3}
\cot^{(1)}_{D \Delta}(z),\ldots,$$ where $(m_0,m_1,m_2,m_3) \in
\N_0^4$, or finite superpositions for example. Rectangular
multi-periodic instanton solutions for the cases $p=2,3,4$ are
obtained by considering the generalized $p$-fold periodic
cotangent functions (see also \cite{LR2}) $$ \cot^{(p)}_{D
\Delta}(z) = z^{-1} + \sum_{\mathbb{Z}^p
\backslash\{0\}}\Big[(z-w)^{-1}+\sum_{\mu=0}^{3}
(w^{-1}z)^{\mu}w\Big] \;\; {\rm for~} p=2,3,4.$$ The application of
the quaternionic residue theorem leads to
\begin{eqnarray*}
\frac{1}{8 \pi^2} \int_{\partial P} d\sigma(z) \cot^{(4)}_D(z) &=&
\frac{1}{8 \pi^2} \int_{\partial P} d\sigma(z) \Delta
\cot^{(4)}_{D \Delta}(z) = \sum_{a \in P} {\rm res}
(\cot^{(4)}_D;a)\\ &=& {\rm res} (\cot^{(4)};0) = 1,
\end{eqnarray*}
where $\cot^{(4)}_D$ denotes the four-fold periodic monogenic
cotangent function discussed in \cite{KraHabil}.
 The choice of the ansatz $$ a_{\mu} = \frac{1}{2}
Vec(e_{\mu} \Delta \cot^{(4)}_{D \Delta}(z) (D \cot_{D
\Delta}(z))^{-1})$$ corresponds to a Yang-Mills solution with unit
instanton number per period parallelepiped. The same result can be
established for $p < 4$ with a similar argument. Since
$\cot^{(p)}_{D \Delta}$ is the fundamental solution of $D \Delta$,
it provides us with all solutions to $D \Delta$ on
$\mathbb{H}/\mathbb{Z}^p$. This can be seen by applying the
following Green's integral formula taken from \cite{KraRyan2}. In
the case $p=1$ we have
\begin{eqnarray*}
  w_{\Gamma'}(y')g(y') &=& \frac{1}{8 \pi^2} \int_{\Gamma'}
\cot^{(1)}_D (z'-y')d\sigma(z')g'(z')\\ & &- \frac{1}{4 \pi^2}
\int_{\Gamma'} \cot^{(1)}_{\Delta}(z'-y') \overline{d\sigma(z')}
D'g'(z')\\ & &  + \frac{1}{2 \pi^2} \int_{\Gamma'} \cot^{(1)}_{D
\Delta}(z'-y') \Delta'g'(z').
\end{eqnarray*}
 In the cases $p=2,3,4$ this formula
can be adapted directly using properly chosen sums of generalized
cotangent functions. By means of the generalized monogenic
cotangent functions we also obtain an argument principle on $C_p$
for isolated zero points. Due to what has been shown in
\cite{KraHabil} Chapter 2.11,
$$
ord(f';c') = \frac{1}{8 \pi^2} \int_{\partial B'(c',\varepsilon)}
\cot^{(p)}_D (f'(z')) \cdot [(Jf)^{ad}(z')]*[d\sigma'(z')].
$$
For $p=4$ the function $\cot^{(4)}$ needs to be slightly modified
like for the Green's integral formula. Here again $``*"$ denotes
the matrix multiplication and $``\cdot"$ the quaternionic
multiplication. The vector
$[(Jf)^{ad}(z')]*[d\sigma'(z')]\in\mathbb{R}^4$ is again
re-interpreted as a quaternion. This formula provides a direct
relation between the fundamental solution and the Chern number of
the $SU(2)$ principle bundles with base manifolds ~$C_p$.

\subsection{The Hopf manifold}

According to \cite{KraRyan2} the fundamental solution of $D'$ and
$D\Delta'$ on $S^3 \times S^1$ with values in the trivial spinor bundle are given by
$$
G^H_D(z',y') = \sum_{k=-\infty}^{0} m^{3k/2}q_{\bf 0}(m^k z' -y')
+ q_{\bf 0}(z') [\sum_{k=1}^{+\infty} m^{-3k/2}q_{\bf
0}(m^{-k}z'^{-1}-y'^{-1})]q_{\bf 0}(y'),
$$
$$
G^H_{D \Delta}(z',y') = \sum_{k=-\infty}^{0} m^{k/2}(m^k z'
-y'^{-1}) + z'^{-1} [\sum_{k=1}^{+\infty}
m^{-k/2}(m^{-k}z'^{-1}-y'^{-1})^{-1}]y'^{-1}.
$$
These induce self-dual $SU(2)$-instanton solutions on $S^3 \times
S^1$. By applying the Green's formula we obtain all solutions of
$D \Delta$ on $S^3 \times S^1$. For this we have to replace the
generalized cotangent functions from the previous formula by the
projection of $G^H$. This yields an argument principle for
isolated zeroes for monogenic functions on $S^3 \times S^1$. More
precisely, we can establish
\begin{theorem}
Let $G'\subseteq S^3 \times S^1 =:H$ be a domain. Let
$f':G'\rightarrow H$ be left monogenic in $G'$ and $c'\in G'$ be
an isolated zero point. Let $\varepsilon > 0$ be a real such that
$f'|_{B'(c',\varepsilon)} \neq 0$, and that
$B'(c',\varepsilon)\subset G'$, where $B'(c',\varepsilon)$ denotes
the projection of the ball $B(c,\varepsilon)$ onto the Hopf
manifold $H$. Then we have
$$
ord(f',c') = \frac{1}{8 \pi^2} \int_{\partial
B(c',\varepsilon)}G^H_D(f'(z'))\cdot
[(Jf)^{ad}(z')]*[d\sigma'(z')].
$$
\end{theorem}

In order to show that ord $(f';c')$ actually is an integer, we
need to use the generalized Cauchy integral formula on the Hopf
manifold. The rest of the argumentation is then similar as in the
Euclidean case treated in \cite{HeKra}. In view of the
transformation rule for differential forms from \cite{Zoll}
reading $d\sigma'(f'(z'))=[(Jf')^{ad}(z')]
* [d\sigma'(z')]$ the integral above can be written in the form
$$
ord(f',c') = \frac{1}{8 \pi^2} \int_{\partial
f(B(c',\varepsilon))}G^H_D(z') d\sigma'(z').
$$
The generalized Cauchy integral formula for left monogenic
functions on the Hopf-manifold tells us that
$$
\frac{1}{8 \pi^2} \int\limits_{\Gamma} G^H_D(z',y') d\sigma(z')
g(z') = w_{\Gamma'}(y') g(y')
$$
Now we replace  $g'\equiv 1$ and $y'$ by $f'(c')$ and $\Gamma'$ by
$f'(\partial B'(c',\varepsilon))$. This leads consequently to
$$
\frac{1}{8 \pi^2}\int\limits_{f'(\partial
B'(c,\varepsilon))}G^H_D(z',\underbrace{f'(c')}_{=0})d\sigma'(z')=w_{f'(\partial
B'(c',\varepsilon))}(\underbrace{f'(c')}_{=0}).
$$
The value ord$(f',c')=w_{f'(\partial B'(c',\varepsilon))}(0)$ is
thus the integer counting how often the image of $B'$ under $f'$
around the isolated zero wraps around zero.

\par\medskip\par

Similarly to the Euclidean case, one can replace the projection of
the ball by a null-homologous $3$-dimensional cycle parameterizing
an $3$-dimensional surface of a $4$-dimensional simply connected
domain inside of $G'\subset H$ which contains the isolated zero
$c'$ in its interior and no further zeroes neither in its interior
nor on its boundary.

\par\medskip\par

 Again, this allows to relate the fundamental solution on $S^3\times S^1$ and the Chern numbers of 
 $SU(2)$ principle bundles.

\subsection{k-handled cylinders and tori}

According to~\cite{BCKR}, the fundamental solution of $D \Delta$
on ${\cal{M}}=H^{+}/\Gamma_p[N]$ ($N \ge 3$) is induced by the Fueter
holomorphic Poincar\'e series on $\Gamma_p[N]$. Adapting the representation formulas of \cite{BCKR} to the particular context considered here, the 
corresponding Poincar\'e series then reads
$$
G_{{\cal{M}}}(z,y) = \frac{1}{8 \pi^2} \sum_{M \in \Gamma_p[N]}
x_n^2 \frac{\overline{cz+d}}{|cz+d|^6}(M<z>-y)^{-1}.$$
Similar to
the previously discussed example, this series provides us with the
Green's kernel functions. Likewise, a residue and argument formula
on this class of conformally flat manifolds can be obtained which
can be used to express the Chern number of $SU(2)$ principal
bundles over k-handled cylinders and tori in terms of the above
given fundamental solution. Notice that we exploited in the proof to show that the topological winding number (which gives 
the order of an $a$ point) is an integer that the constant function $g \equiv 1$ is a solution to the 
Cauchy-Riemann operator on the manifold. This is true since we are dealing with conformally flat mainfolds where the
 scalar curvature is zero. In the case of 
dealing with other manifolds where we have a non-zero scalar 
curvature we have to apply a more sophisticated argument to relate the 
second Chern number with the fundamental solution of $D\Delta$.\\[.3cm] 

{\bf Acknowledgements}: The authors are very thankful to Professor John Ryan from the University of Arkansas and to Professor Vladimir Soucek from Charles University of Prague for the very fruitful discussions which lead to a successful development of this paper. 

The research leading to these results has received funding from the
European Research Council under the European Union's Seventh Framework
Programme (FP7/2007-2013) / ERC grant agreement n$^\circ$~267087.

\end{document}